\documentclass{article}
\usepackage{amsmath,amssymb}

\newtheorem{theorem}{Theorem}
\newtheorem{definition}{Definition}
\newtheorem{lemma}[theorem]{Lemma}

\newtheorem{con}[theorem]{Conjecture}

\newcommand{\inte }{{\rm int}\,}
\newcommand{\diam }{\,{\rm diam}\,}

\def\qed{{\hfill{\vrule height5pt width3pt depth0pt}\medskip}}
\def\comment#1{{}}

\begin{document}
\begin{center}
{\Large \bf  A Lohner-type algorithm for control systems and
ordinary differential inclusions}

 \vskip 0.5cm
{\large Tomasz Kapela}, {\large Piotr Zgliczynski}\footnote{
Research supported in part by Polish State Ministry of
            Science and Information Technology  grant N201 024 31/2163 
} \vskip 0.2cm
  Jagiellonian University, Institute of Computer Science, \\
 Nawojki 11, 30--072  Krak\'ow, Poland \\ e-mail: kapela@ii.uj.edu.pl, zgliczyn@ii.uj.edu.pl

\vskip 0.5cm

\today
\end{center}

\begin{abstract}
  We describe  a  Lohner-type algorithm for the computation of rigorous upper bounds
  for reachable set for   control systems,
  solutions of ordinary differential inclusions
  and perturbations of ODEs.
\end{abstract}


\section{Introduction}

Our goal is to  present a Lohner-type algorithm for an rigorous
integration of perturbations of ODEs, which can be seen also as an
algorithm for an integration of control systems or ordinary
differential inclusions. This paper depends heavily on
\cite{ZC1Lo}, as the proposed algorithm is  a modification running
on top of the $C^0$-Lohner algorithm for ODEs described (after
\cite{Lo,Lo1}) there.

We study the following nonautonomous  ODE
\begin{eqnarray}
  x'(t)=f(x(t),y(t)), \quad x(0) = x_0  \label{eq:full}
\end{eqnarray}
where $x \in {\mathbb R}^n$, $f:\mathbb{R}^n \times \mathbb{R}^m
\to \mathbb{R}^n$ is $C^1$  and $y:\mathbb{R} \supset D \to
{\mathbb R}^m$. Assume that we have some knowledge about $y(t)$,
for example $|y(t)| < \epsilon$ for $0 \leq t \leq T$. We would
like to find an rigorous enclosure for $x(t)$.

The problem of this type arises, for example, in the context of
the control theory (see \cite{G, KS, Sz}) and in the rigorous
integration of dissipative PDEs (see \cite{ZM, ZKSper,ZKS3} for
more details). In this last setting $x$ represents the dominating
modes and $y$ is a tail of the Fourier expansion, so that
(\ref{eq:full}) is complemented by the equation for $y$ of the
form $y'(t)=g(x(t),y(t))$ for which we are able to produce some a
priori bounds. The proposed algorithm works, as we were able using
it prove the existence of multiple periodic orbits for
Kuramoto-Sivashinsky PDE \cite{ZKSper,ZKS3}.

The proposed algorithm can also be used to find  rigorous bounds
for solutions of differential inclusions
\begin{equation}
  x' \in h(x) + \epsilon(t),  \label{eq:diffincl}
\end{equation}
where $h:{\mathbb R}^n \to {\mathbb R}^n$ is a $C^1$-vector field
and $\epsilon(t) \subset {\mathbb R}^n$. We can cast
(\ref{eq:diffincl}) in the form (\ref{eq:full}) by setting
$f(x,y)=h(x)+y$ and requiring that $y(t) \in \epsilon(t)$ for all
$t$.

Non-autonomous ODEs represent another important class  of
applications. While one can easily modify the Lohner algorithm to
handle a non-autonomous ODE directly, it makes sense to apply the
proposed  Lohner-type algorithm for perturbed ODEs  for
(\ref{eq:full}), because only in this way we can estimate
rigorously the Poincar\'e map on a section $\alpha(x)=0$ (defined
in terms of $x$ only) for any initial conditions $(x,t_0)$. This
kind of algorithm shall allow to attack the question of symbolic
dynamics for non-autonomous ODEs (see \cite{CZ}) and ODEs with
small delays (see \cite{WjZ}).

Another new element in this paper, besides the proposed
algorithm, is a new inequality concerning bounds for perturbations
of ODEs. It is contained  in  Theorem~\ref{thm:alg-comp-wise} and
links together the component-wise estimates based on one-sided
Lipschitz conditions (see \cite{W}) and the logarithmic norms
(see \cite{D,L}).

The content of the present paper can be described as follows: in
Section~\ref{sec:control} we define a notion of weak solution of
(\ref{eq:full}) and state some facts from the theory of Lebesgue
integration. In Section~\ref{sec:lognorm} we recall the notion of
the logarithmic norm and state its basic properties. In
Sections~\ref{sec:compestm} and~\ref{sec:formul}  we derive basic
estimates for comparison of perturbed and unperturbed ODEs. In
Section~\ref{sec:alg} we give a description of one step of the
proposed Lohner-type algorithm. In Section~\ref{sec:betwper} we
describe how to estimate the trajectory of (\ref{eq:full}) between
time steps which allows to compute the Poincar\'e map. In the
following section we discuss some tests.

The algorithm presented in this paper was implemented as a part of CAPD library (see \cite{CAPD}).
This library contains many tools for rigorous computations and computer assisted proofs in the contexts of dynamical systems. 
All the tests in Section~\ref{sec:testexamples} was performed using CAPD library.

\subsection{Basic notation } We will use the same conventions as in
\cite{ZC1Lo}. In the sequel, by arabic letters we denote single
valued objects like vectors, real numbers, matrices. Quite often
in this paper we will use square brackets, for example $[r]$, to
denote sets. Usually this will be some set constructed in the
algorithm. Sets will also   be denoted by single letters, for
example $S$, when it is clear from the context that it represents
a set. In situations when we want to stress (for example in the
detailed description of algorithm) that we have  a set in a
formula involving both single-valued objects and sets  we will
rather use the square bracket, hence we prefer to write $[S]$
instead of $S$ to represent a set. From this point of view $[S]$
and $S$ are different symbols in the alphabet used to name
variables and formally speaking there is no relation between the
set represented by $[S]$ and the object represented by $S$. Quite
often in the description of the algorithm we will have a situation
that both variables $[S]$ and $S$ are used simultaneously, then
usually $S \in [S]$, but this is always stated explicitly.

For a set $[S]$ by $[S]_I$ we denote the interval hull of $[S]$,
i.e. the smallest product of intervals containing $[S]$. The
symbol $\mbox{hull}(x_1,\dots,x_k)$ will denote the interval hull
of intervals $x_1,\dots,x_k$. For any interval set $[S]=[S]_I$ by
$\mbox{m}([S])$ we will denote a center point of $[S]_I$. For any
interval $[a,b]$ we define a diameter by $\mbox{diam}([a,b])=b-a$.
For an interval vector or an interval matrix  $[S]=[S]_I$ by
$\diam([S])$ we will denote the maximum of diameters of its
components. For an interval $[x^-,x+]$ we set
$right([x^-,x^+])=x^+$ and $left([x^-,x^+])=x^-$.

If $f(x_1,x_2,\dots,x_k)$ is a function and let
$X_1,X_2,\dots,X_k$ be some sets, then by
\begin{displaymath}
f(X_1,\dots,X_k) = \{ f(z_1,\dots,z_k)\: | \: \mbox{ where $z_i
\in X_i$ for $i=1,\dots,j$}\}
\end{displaymath}

 For a set $X \subset {\mathbb R}^d$ by $\inte X$ we
denote an interior of $X$. For $\mathbb{R}^n$  we will denote the
norm of $x$ by $\|x\|$ and if the formula for the norm is not
specified in some context, then it means that it is ok to use any
norm  there. Let $x_0 \in \mathbb{R}^s$, then $B(x_0,r)=\{z \in
\mathbb{R}^s \: | \: \|x_0 - z \| < r \}$.

For $v , w \in {\mathbb R}^n$ and $A,B \in {\mathbb R}^{n\times
n}$ ($n=1,\dots,\infty$) we say that
\begin{eqnarray*}
  v \leq w \qquad &\mbox{iff}& \qquad  \forall i \quad v_i \leq w_i,   \\
  A \leq B \qquad &\mbox{iff}& \qquad \forall
  ij \quad  A_{ij} \leq B_{ij}.
\end{eqnarray*}


\subsection{Warning.}
At the first encounter with the question of an rigorous
integration of (\ref{eq:full}) one may hope that the direct
application of any algorithm for rigorous integration of ODEs
should be enough for (\ref{eq:diffincl}). To this end consider a
differential inclusion
\begin{equation}
  x' \in f(x) + [\epsilon ],  \qquad [\epsilon] = \Pi_{i=1}^n[-\epsilon_i,\epsilon_i]
  \label{eq:diffincl1}.
\end{equation}
and a related ODE
\begin{equation}
  x'=f(x) + \epsilon, \qquad \epsilon \in [\epsilon].
  \label{eq:odeper}
\end{equation}

 One may naively hope that, for example, the Lohner algorithm applied to
 (\ref{eq:odeper}) with $[\epsilon]$ as an interval parameter
 in the definition of a constant term in $f(x)$ will give an
 enclosure not only for (\ref{eq:odeper}), but also for
 (\ref{eq:diffincl1}). For this to be true we need the
 following
 \begin{con}
   \label{lem:unvalid}
   Assume $x(t)$ satisfies (\ref{eq:diffincl1}) for $t \in [0,T]$.

   Then for any $t \in [0,T]$ there exists $\epsilon \in [\epsilon]$ such that
   $x_\epsilon(t)=x(t)$ and $x_\epsilon(0)=x(0)$, where
   $x_\epsilon$is a solution of (\ref{eq:odeper}).
 \end{con}
The above conjecture is false as shown by the following example
\cite{Se}.

Consider a differential inclusion given by
\begin{eqnarray}
  x'&\in &y + [-\epsilon,\epsilon],  \label{eq:oscy_per} \\
  y' & \in &-x + [-\epsilon,\epsilon]. \nonumber
\end{eqnarray}
For fixed $\delta \in [-\epsilon,\epsilon]^2$ we have the
following system of ODEs
\begin{eqnarray}
  x'&=&y + \delta_1,  \label{eq:oscy} \\
  y' &=&-x + \delta_2, \nonumber
\end{eqnarray}
all solutions with an initial condition in a compact set have a
uniform bound independent of $\delta$ for  $t
>0$, which is given by the energy integral for (\ref{eq:oscy})
\begin{equation}
  (x- \delta_2)^2 + (y + \delta_1)^2.
\end{equation}
This is not the case for the solutions of (\ref{eq:oscy_per}) as
it is clearly seen for $\epsilon(t)$ given as  a resonant forcing
\begin{eqnarray}
  x'&=&y ,  \label{eq:oscy_res} \\
  y' &=&-x + \epsilon \sin t. \nonumber
\end{eqnarray}

\section{Control Systems, the notion of the solution}
\label{sec:control}

In this section we define a notion of (weak) solution of
(\ref{eq:full}).

We use some standard notions from the measure theory, see for
example \cite{Ru} for precise definitions. The integral will
always mean the Lebesgue integral and the measure of the set is
always Lebesgue measure.

\subsection{Some facts from the theory of Lebesgue integral}

We will denote by $m(E)$ the Lebesgue measure of $E$.

Let $D$ be a measurable subset of $\mathbb{R}^k$. By $L^1(D)$ we
will denote a set of measurable functions $f:D \to \mathbb{R}$
such that $\int_D |f| dm < \infty$. If $f:D \to \mathbb{R}^n$ is
measurable, then we say that $f \in L^1(D)$ if function $\| f\|
\in L^1(D)$.

\begin{definition}
Let $D \subset \mathbb{R}$ be  an interval. Function $f:D \to
\mathbb{R}^k$ is \emph{absolutely continuous }, if for every
$\epsilon > 0$ there exists $\delta > 0$, such that for any family
of disjoint intervals
$(\alpha_1,\beta_1),\dots,(\alpha_N,\beta_N)$ such that
\begin{displaymath}
  \sum_{i=1}^N (\beta_i - \alpha_i) < \delta
\end{displaymath}
the following inequality is satisfied
\begin{displaymath}
  \sum_{i=1}^N (f(\beta_i) - f(\alpha_i)) < \epsilon
\end{displaymath}
\end{definition}

The following statement follows directly from results about the
differentiability of measures and functions of bounded variation
(see \cite[Chapter 8]{Ru}).

\begin{theorem}
\label{thm:acdiff}
 Let $D=[a,b]$, $x:D \to \mathbb{R}^n$.

 There exists  $g:D
\to \mathbb{R}^n$ a measurable function such that equation
\begin{equation}
   x(t)-x(a)=\int_a^t g(s)ds
\end{equation}
holds for all $t \in [a,b]$ iff $x$ is absolutely continuous.
 In
this situation $x'(t)$ exists almost everywhere in $[a,b]$ and
$x'(t)=g(t)$.
\end{theorem}

\begin{definition}
Assume $x \in \mathbb{R}^k$. We say that a sequence $\{E_i\}$ of
measurable subsets of $\mathbb{R}^k$ \emph{converges well to the
point $x$}, if there exists $\alpha >0$ such that: \\
every set $E_i$ is contained in $B(x,r_i)$, such that
\begin{equation}
  m(E_i) \geq \alpha m(B(x,r_i)),\quad  \lim_{i \to \infty}r_i=0
\end{equation}

\end{definition}

In the sequel we will need the following theorem
\begin{theorem}
\cite[Thm. 8.8]{Ru} Assume that $f \in L^1(\mathbb{R}^k)$ and
define a Lebesgue set $L_f$ of the function $f$ as the set of all
points $x_0 \in \mathbb{R}^k$ for which
\begin{equation}
  \lim_{i\to \infty} \frac{1}{m(E_i)}\int_{E_i}|f(x) - f(x_0)|
  dx=0
\end{equation}
for every sequence $\{E_i\}$ converging well to the point $x_0$.

Then set $L_f$ contains almost all points of $\mathbb{R}^k$.
\end{theorem}

The above theorem immediately implies the following lemma.
\begin{lemma}
 \label{lem:smallinte}
  Let $f:[a,b] \to \mathbb{R}^k$ be a measurable function. Then
  for almost all points  $ x \in [a,b)$ holds
  \begin{equation}
    \lim_{h \to 0^+}\frac{1}{h} \int_{x}^{x+h} \|f(s) - f(x)
    \|ds=0
  \end{equation}
\end{lemma}

\subsection{Weak solutions of ODEs}

Control System is given by equation
\begin{eqnarray}
  x'(t)=f(x(t),y(t)) \quad x(t_0)=x_0  \label{eq:control_system}
\end{eqnarray}
where $x \in {\mathbb R}^n$, $f:\mathbb{R}^n \times \mathbb{R}^m
\to \mathbb{R}^n$ is $C^1$ and $y:\mathbb{R} \supset D \to
{\mathbb R}^m$ is a measurable function from a given class $U$.

Because the right hand side of (\ref{eq:control_system}) can be
non-continuous we need to define what we mean by solution of
(\ref{eq:control_system}).

\begin{definition}
\label{def:solOde}
 Let $D \subset \mathbb{R} $ be an interval (a connected subset of $\mathbb{R}$) containing $t_0$.

 An absolutely continuous function $x:D \to \mathbb{R}^n$ is a \emph{weak solution}
 of
  (\ref{eq:control_system}) if for all $t \in D$ holds
\begin{equation}
 x(t) = x_0 + \int_{t_0}^t f(x(s), y(s)) ds.  \label{eq:inteq}
\end{equation}

We say that a continuous function $x: D \to \mathbb{R}^n$ is a
\emph{(classical) solution } of (\ref{eq:control_system}) if
$x'(t)$ exists for all $t \in \inte D$, $x(t_0)=x_0$ and
\begin{equation}
  x'(t)=f(x(t),y(t)), \qquad \forall t \in \inte D.
\end{equation}
\end{definition}

From Theorem~\ref{thm:acdiff}  it follows that  $x$ is a weak
solution of (\ref{eq:control_system}) iff
\begin{equation}
  x'(t)=f(x(t),y(t)), \quad \mbox{allmost everywhere in $D$}
\end{equation}
and the function $t \mapsto f(x(t),y(t))$ is in $L^1(D)$. Hence
the weak solution in the sense  of Def.~\ref{def:solOde} is a
solution of (\ref{eq:control_system}) in the sense of Caratheodory
\cite{W}.

In the remainder of this paper we will always consider the
function $f$ on the right hand side of (\ref{eq:control_system})
to be of class $C^k$ (for $k \geq 1$) and $y$ to be  bounded on
compact intervals and  measurable. In such situation the integral
equation (\ref{eq:inteq}) has a unique solution defined for  some
$h >0$ on  $[t_0,t+h]$. The proof of this fact is a
straightforward application of the Banach contraction principle
\cite{W}.

\section{Basic facts on logarithmic norms}
\label{sec:lognorm}

 Let $\| \cdot \|$ denote a vector norm on $\mathbb{R}^n$  as
well as its subordinate matrix (operator) norm on $\mathbb{R}^{n
\times n}$. The classical definition of the \emph{logarithmic
norm} of matrix $A$,
\begin{equation}
  \mu(A)=\lim_{h \to 0^+} \frac{\|I + hA \| - 1}{h} \label{eq:lnormdef}
\end{equation}
was introduced in 1958 independently by Dahlquist \cite{D} and
Lozinskii \cite{L}.

In this section we will briefly recall some basic facts, with
proofs, about the logarithmic norms. For survey regarding the
modern developments stemming from this notion the reader is
referred to \cite{So} and the literature given there. Our
presentation is based on \cite[Ch. 1.5 ]{DV}, which was based on
\cite{D}.

\begin{lemma}
\label{lem:lognormexists}
  For any matrix $A \in \mathbb{R}^{n \times n}$. The limit in (\ref{eq:lnormdef})
  exists and
  \begin{eqnarray}
    \frac{\|I + h_1A \| - 1}{h_1} &\leq& \frac{\|I + h_2A \| -
    1}{h_2}, \qquad \mbox{for $0 < h_1 < h_2$}  \label{eq:lnmonot}\\
     -\|A\| &\leq&  \mu(A) \leq \|A \|.
  \end{eqnarray}
\end{lemma}
\textbf{Proof:} Let us fix $h >0$ and let $0 < \theta <1$, then
\begin{eqnarray*}
  \|I + \theta h A\| = \|\theta(I + hA) + (1-\theta)I\| \leq
  \theta \|I + hA\| + (1-\theta)\|I\|.
\end{eqnarray*}
From this immediately obtain
\begin{equation}
  \frac{\|I + \theta h A\| -1}{\theta h} \leq   \frac{\|I + h A\|
  -1}{h},
\end{equation}
which proves (\ref{eq:lnmonot}).

From the triangle inequality one gets
\begin{equation}
  -h \| A \| \leq \|I + h A \| - \| I \| \leq h\| A \|,
\end{equation}
therefore
\begin{equation}
 -\| A\|  \leq   \frac{\|I
+ h A\|  -1}{h} \leq \| A \|.
\end{equation}
The monotonicity (\ref{eq:lnmonot}) and the existence of the lower
bound imply the existence of $\mu(A)$.
 \qed

\begin{theorem}
\label{thm:lognormcont}
  The function  $\mu: \mathbb{R}^{n \times n} \to \mathbb{R}$, which assigns to $A$
  its logarithmic norm is continuous and convex. Moreover,  functions $\mu(h,A)=\frac{\|I + hA  \| -1}{h}$
  converge locally uniformly and monotonically to $\mu(A)$  for $h \to
  0^+$.

  To be more precise, for any compact set $K \subset \mathbb{R}^{n \times
  n}$ and any $\epsilon >0$ there exists $h_0 >0$, such that for
  all $0< h < h_0$ and any $A \in K$ holds
  \begin{equation}
      \epsilon > \mu(h,A) - \mu(A) \geq 0.
  \end{equation}
\end{theorem}
\textbf{Proof:} Let $h >0$. An easy computation show that,  for
any $0 \leq \lambda \leq 1$ and $A_1,A_2 \in \mathbb{R}^{n\times
n}$ holds
\begin{eqnarray*}
  \mu(h, \lambda A_1 + (1-\lambda) A_2) \leq  \lambda \mu(h,A_1) +
  (1-\lambda) \mu(h,A_2).
\end{eqnarray*}
Therefore, for any $h>0$ function $\mu(h,\cdot):\mathbb{R}^{n
\times n} \to \mathbb{R}$ is convex.

By taking the limit $h \to 0^+$ from Lemma~\ref{lem:lognormexists}
it follows that $\mu(A)$ is a convex function. Observe that on any
bounded set $U \subset \mathbb{R}^{n \times n}$   $\mu(A)$ is
bounded by $\sup_{A \in U} \|A\| < +\infty$, therefore from the
theory of convex functions (see for example \cite[Chap. 6]{La}) it
follows that $\mu$ is continuous. The uniform convergence of
$\mu(h,\cdot)$ to $\mu$ on compact sets follows from Dini's
Theorem on monotone sequences of pointwise converging continuous
functions to continuous limit and Lemma~\ref{lem:lognormexists}.
\qed

The following lemma follows directly from the convexity of
$\mu(A)$
\begin{lemma}
\label{lem:lgnornconcave} Let $A:[0,1] \to \mathbb{R}^{n \times
n}$ be a bounded measurable function. Then
\begin{equation}
  \mu\left(\int_{0}^1 A(s) ds\right) \leq \int_0^1 \mu(A(s)) ds
  \leq \sup_{s \in [0,1]}
  \mu(A(s)).
\end{equation}
\end{lemma}

\section{Bounds for perturbations of  ODEs}
\label{sec:compestm}

In this section we state the basic theorem comparing a solution of
an ODE and an approximate solution. Our approach unifies the
approach based on logarithmic norms and one-sided Lipschitz
condition leading to component-wise bounds from \cite[Ch.
II.13]{W}.

\subsection{Estimates for non-autonomous linear equations}
Consider a linear equation
\begin{equation}
  x'(t)=A(t) \cdot x(t) + b(t), \label{eq:linnonhp}
\end{equation}
where $x(t) \in \mathbb{R}^k$, $A(t) \in \mathbb{R}^{k \times k}$,
$b(t) \in \mathbb{R}^k$, $A$ and $b$ are bounded and measurable.

We would like give some bounds on solutions of
(\ref{eq:linnonhp}). We assume that our phase space $\mathbb{R}^k$
is decomposed as follows $\mathbb{R}^k=\oplus_{i=1}^n
\mathbb{R}^{k_i}$. Therefore,  we have a decomposition of $z \in
\mathbb{R}^k$ into $(z_1,\dots,z_n)$ such that $z_i \in
\mathbb{R}^{k_i}$. In this section we will carefully distinguish
between the symbol $\| \cdot \|$ and $|\cdot|$. The symbol $\|
\cdot \|$ will always denote a norm, but the symbol $|z|$ for $z
\in \mathbb{R}^k$  will usually denote a vector of norms of $z_i$,
but this will be always clearly indicated in the text. Observe
that, when we have such decomposition, then equation
(\ref{eq:linnonhp}) can be written as follows
\begin{equation}
  z'_i(t)= \sum_{j} A_{ij}(t)z_j(t) + b_i(t), \quad i=1,\dots,n
\end{equation}
where $z_i, b_i \in \mathbb{R}^{k_i}$ and $A_{ij}(t) \in
L(\mathbb{R}^{k_i},\mathbb{R}^{k_j})$ is a linear map (a matrix).
In this way matrix $A$ is decomposed into blocks $A_{ij}$. For
each block we will assign number $J_{ij}$ and collect them in
matrix $J$. Roughly speaking $J_{ij}$ will estimate an influence
of $z_j$ on $z'_i$.

The fundamental lemma in this section is:
\begin{lemma}
\label{lem:comptestm} Assume that $z:[0,T] \to {\mathbb
R}^k=\oplus_{i=1}^n \mathbb{R}^{k_i}$ is an absolutely continuous
map, which is a weak solution of the equation
\begin{equation}
  z'(t)=A(t) \cdot z(t) + \delta(t),
\end{equation}
where $\delta:[0,T] \to {\mathbb R}^k$ and $A:[0,T]\to {\mathbb
R}^{k \times k}$ are bounded and measurable.

Assume that  measurable  matrix function  $J:[0,T] \to
\mathbb{R}^{n \times n}$ satisfies the following inequalities for
all $t \in [0,T]$
\begin{equation}
\label{eq:def-J}  J_{ij}(t) \geq
  \begin{cases}
      \|A_{ij}(t) \|  & \text{for $i \neq j$}, \\
     \mu(A_{ii}(t))  & \text{for $i=j$}.
  \end{cases}
\end{equation}

Let $C_i(t) =\|\delta_i(t)\|$ and
$|z|(t)=(\|z_1(t)\|,\|z_2(t)\|,\dots,\|z_n(t)\|)$.

Then
\begin{equation}
  |z|(t) \leq y(t)
\end{equation}
where  $y:[0,T] \to \mathbb{R}^n$ is a weak solution of the
problem
\begin{equation}
   y'(t)=J(t)y(t) + C(t), \qquad y(0)=|z|(0).
\end{equation}
\end{lemma}
\textbf{Proof:} Observe that for all $i$ the function $t \mapsto
\|z_i(t)\|$ is absolutely continuous. Therefore from Theorem
\ref{thm:acdiff} it follows that for almost every $t \in [0,T]$
the derivative of $\|z_i\|$ exists. We will estimate this
derivative for such $t$.

We have
\begin{eqnarray*}
  z(t+h) = z(t)+\int_t^{t+h}A(s)z(s) ds + \int_t^{t+h} \delta(s)
  ds = \\  z(t) + h\left(A(t)z(t)) + h\delta(t)\right) +
   \int_t^{t+h}\left(A(s)z(s) - A(t)z(t)\right)
       +  (\delta(s) - \delta(t)) ds
\end{eqnarray*}

Let us fix $i$ and $t\in [0,T)$. We consider the projection onto
$i$-th subspace. We have
\begin{eqnarray*}
 \|z_i(t+h)\| \leq \|I + h A_{ii}(t)\| \cdot \|z_i\|(t) + h \sum_{j \neq
 i}\|A_{ij}(t)\|\cdot \|z_j(t)\| + h \|\delta_i(t)\| + \\
   \int_t^{t+h}\left\|A(s)z(s) - A(t)z(t)\right\| ds
       +   \int_t^{t+h}\left\|\delta(s) - \delta(t)\right\| ds
\end{eqnarray*}
and then we obtain for $h >0$
\begin{eqnarray*}
 \frac{\|z_i(t+h)\| - \|z_i(t)\|}{h}\leq \frac{\|I + h A_{ii}(t)\|-1 }{h} \cdot \|z_i\|(t)
 + \\ \sum_{j \neq i}\|A_{ij}(t)\|\cdot \|z_j(t)\| +  C_i +
   \frac{1}{h}\int_t^{t+h}\left\|A(s)z(s) - A(t)z(t)\right\| ds
       + \\  \frac{1}{h} \int_t^{t+h}\left\|\delta(s) - \delta(t)\right\| ds
\end{eqnarray*}
Observe that from Lemma~\ref{lem:smallinte} it follows that the
last two terms in the above inequality tend to $0$ as $h\to 0$ for
almost all points in $[0,T)$. From now on we assume that $t$ is
such point.

By passing to the limit with $h \to 0^+$ we obtain for almost all
points in $t \in [0,T]$
\begin{eqnarray}
  \frac{d \|z_i\|}{dt}(t) \leq \mu(A_{ii}(t)) \|z_i\|(t) + \sum_{j \neq
  i}\|A_{ij}(t)\| \cdot \|z_j(t)\| + C_i(t) \leq  \nonumber \\
  \sum_{j}J_{ij}(t) \|z_j\|(t) + C_i(t)
   \label{eq:d|z|}
\end{eqnarray}

Let us define
\begin{eqnarray*}
  x(t)=(x_1(t),x_2(t),\dots,x_n(t))=(\|z_1(t)\|,\|z_2(t)\|,\dots,\|z_n(t)\|),
\end{eqnarray*}

Inequality (\ref{eq:d|z|}) can be rewritten in vector form as
follows
\begin{equation}
  x'(t) \leq J(t) \cdot x(t) + C(t), \label{eq:Dx} \quad \mbox{for almost all $t \in
  [0,T]$}.
\end{equation}
Let  $y:[0,T] \to \mathbb{R}^n$  be a weak solution of
\begin{equation}
  y'(t) = J(t) \cdot y(t) + C(t),  \label{eq:ysol}
\end{equation}
such that $y(0) > |z|(0)=x(0)$.

We want to show that
\begin{equation}
  x(t) < y(t), \quad t \in [0,T]. \label{eq:x-less-y}
\end{equation}

Let us take  diagonal matrix $\Lambda \in \mathbb{R}^{n \times
n}$, such that $\Lambda_{ii} + J_{ii}(t) \geq 0$ for all
$i=1,\dots,n$ and $t \in [0,T]$. Let  us define matrix-valued
function $B:[0,T] \to \mathbb{R}^{n \times n}$  by
\begin{equation}
\label{eq:def-B}  B(t) = \Lambda + J(t).
\end{equation}
Obviously $B_{ij}(t) \geq 0$ for all $t \in [0,T]$.

For any $i=1,\dots,n$ from  (\ref{eq:Dx})  we obtain for almost
all $t \in [0,T]$
\begin{eqnarray}
   x'_i(t) + \Lambda_{ii} x_i(t) \leq  \sum_{j} B_{ij}(t)x_j(t) +
  C_i(t),
\end{eqnarray}
hence
\begin{eqnarray*}
   \frac{d}{dt}\left( e^{\Lambda_{ii} t} x_i(t) \right)\leq e^{\Lambda_{ii} t}
  \left( \sum_{j} B_{ij}x_j(t) + C_i(t) \right).
\end{eqnarray*}
The last inequality has the following vector form
\begin{equation}
  \frac{d}{dt} \left(e^{\Lambda t} x(t)\right) \leq e^{\Lambda t} B(t) x(t) + e^{\Lambda t}
  C(t). \label{eq:d+ex}
\end{equation}
From the above inequality and from Theorem~\ref{thm:acdiff} it
follows that
\begin{eqnarray*}
  e^{\Lambda t} x(t) = e^{\Lambda \cdot 0} x(0) + \int_0^t \frac{d}{dt} \left(e^{\Lambda t}
  x(t)\right)(s) ds \leq x(0) + \\
   \int_0^t e^{\Lambda s} B(s) x(s) + e^{\Lambda s}
  C(s) ds.
\end{eqnarray*}
Hence we obtain
\begin{eqnarray}
  x(t) \leq e^{-\Lambda t} x(0) + \int_0^t e^{-\Lambda (t-s)} \left(
  B(s)  x(s) + C(s)  \right) \: ds \qquad \mbox{for $ t \in [0,T]$} \label{eq:integral_ineq}
\end{eqnarray}

An analogous computation applied to (\ref{eq:ysol}) shows that $y$
satisfies the following integral equation
\begin{equation}
 y(t) =  e^{-\Lambda t}y(0) + \int_{0}^t e^{-\Lambda(t-s)}
 \left(B(s)y(s)+ C(s)  \right)\: ds. \label{eq:integral_eq}
\end{equation}

Now we are ready to prove (\ref{eq:x-less-y}). Let
\begin{equation}
  t_0=\sup\{ t \in [0,T] \: | \: y(s) > x(s), \quad s \in [0,t)\}.
\end{equation}
Obviously from the continuity of $y(t) - x(t)$ it follows that
$t_0>0$.
 From (\ref{eq:integral_eq}) and (\ref{eq:integral_ineq}) we obtain
\begin{displaymath}
  y(t_0) - x(t_0) \geq  e^{-\Lambda t_0}(y(0) - x(0)) + \int_0^{t_0}
  e^{-\Lambda  (t_0-s) } B(s)(y(s) - x(s)) \: ds > 0.
\end{displaymath}
By the continuity  inequality $y(t) > x(t)$ will hold for $t \in
[t_0, t_0 + \epsilon)$ for some $\epsilon >0$. Therefore $t_0=T$.

Hence condition (\ref{eq:x-less-y}) holds. By passing to the limit
$y(0) \to x(0)$ we obtain our assertion.
 \qed

\begin{theorem}
\label{thm:alg-comp-wise} Let $ h>0$. Assume that $f:\mathbb{R}^n
\times \mathbb{R}^m  \to \mathbb{R}^n$ be $C^1$ and $y:[t_0,t_0+h]
\to \mathbb{R}^m$ is bounded and measurable.

Let $[W_y] \subset {\mathbb R}^m$ be convex and such that,
$y([t_0, t_0+h]) \subset [W_y]$.

Let $y_c \in [W_y]$. Assume that $x_1,x_2:[t_0,t_0+h] \to
\mathbb{R}^n$, both absolutely continuous, are  weak solutions of
the following problems, respectively
\begin{eqnarray}
  x_1'&=&f(x_1,y_c), \quad x_1(t_0)=x_0,   \label{eq:iprj}\\
  x_2'&=&f(x_2,y(t)), \quad x_2(t_0)=\bar{x}_0.  \label{eq:inonauto}
\end{eqnarray}

Let $[W_1] \subset [W_2] \subset \mathbb{R}^n$ be convex and
compact and such that
\begin{eqnarray*}
x_1(t) \in [W_1], \quad \ x_2(t) \in [W_2] , \qquad \mbox{for $t
\in [t_0,t_0+h] $}.
\end{eqnarray*}

Then the following inequality holds for $t \in [t_0, t_0 +h]$ and
$i=1,\dots,n$
\begin{equation}
  |x_{1,i}(t) - x_{2,i}(t)| \leq   \left(e^{J(t-t_0)} \cdot (x_0 - \bar{x}_0)\right)_i + \left(\int_{t_0}^t e^{J(t-s)} C \:
  ds \right)_i,
\end{equation}
where
\begin{eqnarray*}
{} [\delta] &=& \{f(x,y_c) - f(x,y) \: | \: x \in [W_1], y
\in [W_y] \}, \\
  C_i &\geq&  \sup \left| [\delta_i] \right| , \quad i=1,\dots,n \\
  J_{ij} &\geq&
  \begin{cases}
     \sup  \mu(\frac{\partial f_i}{\partial x_j}([W_2],[W_y])) & \text{if $i=j$}, \\
    \sup \left\| \frac{\partial f_i}{\partial x_j}([W_2],[W_y])\right\| & \text{if $i \neq j$}.
  \end{cases}
\end{eqnarray*}
\end{theorem}
{\bf Proof:} Let $z(t)=x_1(t)-x_2(t)$. We have for $t \in
[t_0,t_0+h]$
\begin{eqnarray*}
  z(t)=\left(x_1(t_0) + \int_{t_0}^t f(x_1(s),y_c)ds \right) -
  \left( x_2(t_0) + \int_{t_0}^t f(x_2(s),y(s))ds \right) =\\
  z(t_0) + \int_{t_0}^t  \left(f(x_1(s),y_c) -  f(x_2(s),y(s))\right)
  ds.
\end{eqnarray*}
Now observe that
\begin{eqnarray*}
f(x_1(t),y_c) - f(x_2(t),y(t))=
  f(x_1(t),y_c) - f(x_1(t),y(t)) + \\ f(x_1(t),y(t)) -
  f(x_2(t),y(t)) =
   \delta(t) + A(t) \cdot (x_1(t) - x_2(t)),
\end{eqnarray*}
where $\delta (t) \in [\delta]$ is bounded and measurable  and
\begin{displaymath}
A_{ij}(t) = \int_0^1 \frac{\partial f_i}{\partial x_j}\left(x_2(t)
+ s(x_1(t) -x_2(t)) ,y(t)\right) ds
\end{displaymath}
is bounded and measurable matrix.

We obtain
\begin{equation}
  z(t)=z(t_0) + \int_{t_0}^t \left( A(s)z(s) + \delta(s) \right)ds
\end{equation}

To apply Lemma~\ref{lem:comptestm} to the function $z=x_1 - x_2$
to obtain (\ref{eq:vectestm}) we need to show that
\begin{equation}
  J_{ij} \geq
  \begin{cases}
      \sup_{t \in [t_0, t_0+h]}\|A_{ij}(t) \|  & \text{for $i \neq j$}, \\
     \sup_{t \in [t_0, t_0+h]}\mu(A_{ii}(t))  & \text{for $i=j$}.
  \end{cases}
\end{equation}
For the off-diagonal terms we have
\begin{eqnarray*}
 \|A_{ij}(t)\| \leq  \int_0^1 \left\|\frac{\partial f_i}{\partial x_j}\left(x_2(t)
+ s(x_1(t) -x_2(t)) ,y(t)\right)\right\| ds \leq \\
 \sup_{x \in
[W_2], y \in [W_y]} \left\|\frac{\partial f_i}{\partial x_j}(x
,y)\right\| \leq J_{ij}.
\end{eqnarray*}
For the diagonal case we use Lemma~\ref{lem:lgnornconcave}.

The result now follows from Lemma~\ref{lem:comptestm}.
\qed

It is possible to  organize the error estimates slightly
differently, namely estimate $[\delta]$ on $[W_2] \times [W_y]$
instead of on $[W_1] \times [W_y]$, which will produce larger
$[\delta]$, but in the same time  estimate $J$ on $[W_2] \times
\{y_c\}$ instead of $[W_2] \times [W_y]$, which should result in
better $J$, to obtain the following variant of the above theorem.
\begin{theorem}
\label{thm:alg-comp-wise2} The same assumptions and notations as
in Theorem~\ref{thm:alg-comp-wise}.

Then the following inequality holds for $t \in [t_0, t_0 +h]$ and
$i=1,\dots,n$
\begin{equation}
  |x_{1,i}(t) - x_{2,i}(t)| \leq   \left(e^{J(t-t_0)} \cdot (x_0 - \bar{x}_0)\right)_i + \left(\int_{t_0}^t e^{J(t-s)} C \:
  ds \right)_i,
\end{equation}
where
\begin{eqnarray*}
{} [\delta] &=& \{f(x,y_c) - f(x,y) \: | \: x \in [W_2], y
\in [W_y] \}, \\
  C_i &\geq&  \sup \left| [\delta_i] \right| , \quad i=1,\dots,n \\
  J_{ij} &\geq&
  \begin{cases}
     \sup  \mu(\frac{\partial f_i}{\partial x_j}([W_2],y_c)) & \text{if $i=j$}, \\
    \sup \left\| \frac{\partial f_i}{\partial x_j}([W_2],y_c)\right\| & \text{if $i \neq j$}.
  \end{cases}
\end{eqnarray*}
\end{theorem}
\textbf{Proof:} We proceed as in the proof of
Theorem~\ref{thm:alg-comp-wise}. But the difference between
$f(x_1(t),y_c)$ and $f(x_2(t),y(t))$ is computed differently.
Namely,
\begin{eqnarray*}
f(x_1(t),y_c) - f(x_2(t),y(t))=
  f(x_1(t),y_c) - f(x_2(t),y_c) + \\ f(x_2(t),y_c) -
  f(x_2(t),y(t)) =
    A(t) \cdot (x_1(t) - x_2(t)) +  \delta(t),
\end{eqnarray*}
where $\delta (t) \in [\delta]$   and
\begin{displaymath}
A_{ij}(t) = \int_0^1 \frac{\partial f_i}{\partial x_j}\left(x_2(t)
+ s(x_1(t) -x_2(t)) ,y_c\right) ds.
\end{displaymath}
We continue as in the proof of Theorem~\ref{thm:alg-comp-wise}.
\qed


\section{Formulas for various cases}
\label{sec:formul}

In this section we rewrite Theorems~\ref{thm:alg-comp-wise}
and~\ref{thm:alg-comp-wise2} in the form, which will be later used
in our algorithm for the integration of differential inclusions.

\subsection{The estimation of perturbations of ODEs based on logarithmic norms}
From Theorem~\ref{thm:alg-comp-wise2} using the trivial
decomposition consisting of the whole space we obtain the
following lemma.
\begin{lemma}
\label{lem:alg-lognorm} Let $ h>0$. Assume that $f:\mathbb{R}^n
\times \mathbb{R}^m  \to \mathbb{R}^n$ be $C^1$ and $y:[t_0,t_0+h]
\to \mathbb{R}^m$  be bounded and measurable.

Let $[W_y] \subset {\mathbb R}^m$ be convex and such that,
$y([t_0, t_0+h]) \subset [W_y]$.

Let $y_c \in [W_y]$. Assume that $x_1,x_2:[t_0,t_0+h] \to
\mathbb{R}^n$  both absolutely continuous, are  weak solutions of
the following problems, respectively
\begin{eqnarray}
  x_1'&=&f(x_1,y_c), \quad x_1(t_0)=x_0,   \label{eq:iprj2}\\
  x_2'&=&f(x_2,y(t)), \quad x_2(t_0)=\bar{x}_0.  \label{eq:inonauto2}
\end{eqnarray}

Let $[W_1] \subset [W_2] \subset \mathbb{R}^n$ be convex and
compact and such that
\begin{eqnarray*}
x_1(t) \in [W_1], \quad \ x_2(t) \in [W_2] , \qquad \mbox{for $s
\in [t_0,t_0+h] $}.
\end{eqnarray*}

Then for any $t \in [0,h]$ holds
\begin{eqnarray*}
  \| x_2(t_0+t) -x_1(t_0+t) \| \leq \\  \exp(lt)\|x_1(t_0) - x_2(t_0) \|  +
   \exp(l t) \int_{t_0}^{t_0+t} \exp(-ls) \| [\delta] \| ds= \\
     \exp(lt)\|x_1(t_0) - x_2(t_0) \|  +  \frac{\| [\delta] \|}{l}(\exp(lt) - 1)
\end{eqnarray*}
where $l={\rm sup} \left( \mu(\frac{\partial f}{\partial x}
([W_2],y_c) )\right)$, and $\mu$ is the logarithmic norm of the
matrix (see \cite{HNW} for the definition) and
\begin{displaymath}
 {} [\delta] = \{f(x,y_c) - f(x,y) \: | \: x \in [W_2], y \in
[W_y] \}.
\end{displaymath}
\end{lemma}

\subsection{A component-wise estimate }
\label{subsec:comp-wise}

From Theorem~\ref{thm:alg-comp-wise} using the trivial
decomposition $\mathbb{R}^m=\bigoplus_{i=1}^m \mathbb{R}$ we
obtain the following lemma.
\begin{lemma}
\label{lem:alg-comp-wise} Let $ h>0$. Assume that $f:\mathbb{R}^n
\times \mathbb{R}^m  \to \mathbb{R}^n$ be $C^1$ and $y:[t_0,t_0+h]
\to \mathbb{R}^m$ is bounded and measurable.

Let $[W_y] \subset {\mathbb R}^m$ be convex and such that,
$y([t_0, t_0+h]) \subset [W_y]$.

Let $y_c \in [W_y]$. Assume that $x_1,x_2:[t_0,t_0+h] \to
\mathbb{R}^n$, both absolutely continuous, are  weak solutions of
the following problems, respectively
\begin{eqnarray}
  x_1'&=&f(x_1,y_c), \quad x_1(t_0)=x_0,  \\
  x_2'&=&f(x_2,y(t)), \quad x_2(t_0)=\bar{x}_0.
\end{eqnarray}

Let $[W_1] \subset [W_2] \subset \mathbb{R}^n$ be convex and
compact and such that
\begin{eqnarray*}
x_1(t) \in [W_1], \quad \ x_2(t) \in [W_2] , \qquad \mbox{for $s
\in [t_0,t_0+h] $}.
\end{eqnarray*}

Then the following inequality holds for $t \in [t_0, t_0 +h]$ and
$i=1,\dots,n$
\begin{equation}
  |x_{1,i}(t) - x_{2,i}(t)| \leq   \left(e^{Jt} \cdot (x_0 - \bar{x}_0)\right)_i + \left(\int_{t_0}^t e^{J(t-s)} C \:
  ds \right)_i, \label{eq:vectestm}
\end{equation}
where
\begin{eqnarray*}
{} [\delta] &=& \{f(x,y_c) - f(x,y) \: | \: x \in [W_1], y
\in [W_y] \}, \\
  C_i &\geq&  \sup \left| [\delta_i] \right| , \quad i=1,\dots,n \\
  J_{ij} &\geq&
  \begin{cases}
     \sup \frac{\partial f_i}{\partial x_j}([W_2],[W_y]) & \text{if $i=j$}, \\
    \sup \left| \frac{\partial f_i}{\partial x_j}([W_2],[W_y])\right| & \text{if $i \neq j$}.
  \end{cases}
\end{eqnarray*}
\end{lemma}

\section{The Lohner-type algorithm for perturbations of ODEs}
\label{sec:alg}

For a given measurable  and bounded on compact intervals function
$y:[0,\infty) \to \mathbb{R}^m$ let $\varphi(t,x_0,y)$ denotes a
weak solution of equation (\ref{eq:full}) with  initial condition
$x(0)=x_0$. For a given $y_0 \in \mathbb{R}^m$ let ${\overline
\varphi}(t,x_0,y_0)$ be a solution of the following Cauchy problem
\begin{equation}
   x'=f(x,y_0), \quad x(0)=x_0  \label{eq:yfixed}
\end{equation}
with the same initial condition  $x(0)=x_0$. Observe that system
(\ref{eq:yfixed}) is a particular case of (\ref{eq:full}) with
$y(t)=y_0$.

Let $U$ be a some family of functions $y:[0,\infty) \to
\mathbb{R}^m$ which are measurable and are uniformly bounded on
any compact interval, i.e. for any $T>0$ there exists  $M(T)$,
such that for every $y \in U$ and every $t \in [0,T]$ holds
$\|y(t)\| \leq M(T)$.

We are interested in finding rigorous bounds for
$\phi(t,[x_0],[y_0])$, where $[x_0] \subset \mathbb{R}^n$ and
$[y_0] \subset U$. The set $[y_0]$ might be defined as some
dynamical process, in this case we may need to compute something
for each time step, or it can be just given by the specifying the
bounds, for example $ y \in [y_0]$ iff $y(t) \in
[-\epsilon,\epsilon]^m$ and $y$ is measurable.

Below we propose a modification of the original Lohner algorithm
\cite{Lo, Lo1} to treat problem (\ref{eq:full}). Our presentation
 follows the description of the $C^0$-Lohner algorithm presented
in \cite{ZC1Lo}.

\subsection{One step of the algorithm}

In the description below the objects with an index $k$ refer to
the current values and those with an index $k+1$ are the values
after the next time step.

We define
\begin{eqnarray*}
[y_k]=\{ y \in U \: | \:  y(t)=z(t_k+t) \quad \mbox{for some $z
\in [y_0]$} \}.
\end{eqnarray*}
For given $[y] \subset U$ we will also use the following notation
\begin{displaymath}
  [y]([t_1,t_2])=\{ z(t) \: |  \:  z \in [y], t \in [t_1,t_2] \}.
\end{displaymath}

One step of the Lohner algorithm is a shift along the trajectory
of system (\ref{eq:full}) with following input and output data:
\newline {\bf Input data:}
\begin{itemize}
\item $t_k$ is a current time
\item $h_k$ is a time step
\item $[x_k] \subset {\mathbb R}^n $,
    such that $ \varphi(t_k,[x_0],[y_0]) \subset [x_k]$
\item eventually some bounds for $[y_k]$
\end{itemize}

\noindent {\bf Output data:}
\begin{itemize}
\item $t_{k+1}=t_{k} + h_k$ is a new current time
\item $[x_{k+1}] \subset {\mathbb R}^n $,  such that $\varphi(t_{k+1},[x_0],[y_0])
      \subset [x_{k+1}]$
\item eventually some bounds for $[y_0][0,t_{k+1})$.
\end{itemize}

We do not specify here a form (a representation) of sets $[x_k]$.
They can be interval sets, balls, doubletons etc. (see \cite{MZ,
ZC1Lo}). This issue is  very important in handling of the wrapping
effect and is discussed in detail in \cite{Lo,Lo1} (see also
Section 3 in \cite{ZC1Lo}).

One step of the algorithm consists from the following parts:
\begin{description}
\item[1.] \textbf{Generation of  a priori bounds for $\varphi$ and $[y_0]([t_k,t_{k+1}])$}.

 We find a convex and
compact set $[W_2] \subset {\mathbb R}^n$ and a convex set $[W_y]
\subset {\mathbb R}^m$, such that
\begin{eqnarray}
  \varphi([0,h_k],[x_k],[y_k]) \subset [W_2]  \\
 {} [y_k]([0,h_k]) \subset [W_y]
\end{eqnarray}
\item[2.] We fix $y_c \in [W_y]$.
\item[3.] \textbf{Computation of an unperturbed $x$-projection.} We apply one step
  of the $C^0$-Lohner algorithm to (\ref{eq:yfixed}) with a time step
  $h_k$ and an initial condition given by $[x_k]$  and $y_0=y_c$.
  As a result we obtain $[{\overline x}_{k+1}] \subset {\mathbb R}^n$ and a convex and compact set
   $[W_1] \subset {\mathbb R}^n$, such that
  \begin{eqnarray*}
     \overline{\varphi}(h_k,[x_k],y_c) &\subset&    [\overline{x}_{k+1}]\\
     \overline{\varphi}([0,h_k],[x_k],y_c) &\subset& [W_1]
  \end{eqnarray*}
\item[4.] \textbf{Computation of the influence of the perturbation.}
Using formulas from Lemmas~\ref{lem:alg-comp-wise} or
\ref{lem:alg-lognorm} we find a set $[\Delta] \subset {\mathbb
R}^n$, such that
\begin{equation}
   \varphi(t_{k+1},[x_0],[y_0]) \subset  \overline{\varphi}(h_k,[x_k],y_c) + [\Delta].
\end{equation}
Hence
\begin{equation}
  \varphi(t_{k+1},[x_0],[y_0]) \subset
 [x_{k+1}]=[\overline{x}_{k+1}]+ [\Delta]  \label{eq:part4}
\end{equation}
\item[5.] Eventually we do some computation to obtain $[y_{k+1}]$
\end{description}

\subsection{Part 1 -  comments}
In the context of an nonautonomous ODE with small and uniformly
bounded $[\delta]$ we can set $[W_y]={\mathbb R}$. To obtain
$[W_2]$ any rough enclosure procedure devised for ODEs should
work. In the context of a dissipative PDE the whole story is  more
complicated and we refer the interested reader to  \cite{ZKSper}.

\subsection{Part 4 - details}
In Lemmas~\ref{lem:alg-lognorm}  and~\ref{lem:alg-comp-wise}  we
have presented two ways to compute $[\Delta]=[\Delta](h)$ for $0
\leq h \leq h_k$.

\noindent \textbf{An approach based on component-wise estimates}
\begin{itemize}
\item[1.] We set
\begin{eqnarray*}
{}[\delta]&=&[\{f(x,y_c) - f(x,y) \: | \: x \in [W_1], y \in [W_y]
\}]_I \\
  C_i &=&  \mbox{right}(\left| [\delta_i] \right|) , \quad i=1,\dots,n \\
  J_{ij}&=&
  \begin{cases}
     \mbox{right} \left( \frac{\partial f_i}{\partial
     x_i}([W_2],[W_y])\right) & \text{if  $i=j$}, \\
     \mbox{right} \left( \left| \frac{\partial f_i}{\partial x_j}([W_2],[W_y]) \right| \right).
         & \text{if $i \neq j$}.
  \end{cases}
\end{eqnarray*}
\item[2.]$  D = \int_0^h e^{J(h-s)} C \: ds $
\item[3.] $[\Delta_i]=[-D_i,D_i]$, for $i=1,\dots,n$
\end{itemize}

It remains to explain how we compute $\int_0^t e^{A(t-s)} C \:
ds$. First observe that
\begin{equation}
\int_0^t e^{A(t-s)} C \:ds = t \left(\sum_{n=0}^\infty
\frac{(At)^n}{(n+1)!}\right) \cdot C.
\end{equation}
We fix any norm $\|\cdot\|$, such that for any matrix $A = (a_{ij})$ we have $|a_{ij}| \leq \|A\|$.
It is not true for general norm, for example if we take vector norm on $\mathbb{R}^2$ defined by  $\|(x_1,x_2)\| = max\{ \frac{1}{100}x_1, x_2\}$ then associated matrix norm of a matrix $\begin{pmatrix} 0 &100 \\ 0 & 0 \end{pmatrix}$ is equal to 1.
We take for example $L^\infty$-norm, i.e.
$\|x\|_\infty=\max_i |x_i|$ ({\em we should rather chose a norm
for which $\| At \|$ is the smallest one}). Let us set
\begin{eqnarray*}
  {\tilde A}= At, \qquad
  A_m = \frac{{\tilde A}^m}{(m+1)!}.
\end{eqnarray*}
In this notation
\begin{eqnarray*}
  \sum_{m=0}^\infty \frac{(At)^m}{(m+1)!} = \sum_{m=0}^\infty A_m
  \\
  A_0 = \mbox{Id}, \qquad A_{m+1}=A_m \cdot \frac{{\tilde A}}{m+2}
\end{eqnarray*}
For the remainder term we will use the following estimate
\begin{eqnarray*}
  \|A_{N + k}\| \leq \|A_N\|  \cdot \left\| \frac{{\tilde A}}{N+2}\right
  \|^k
\end{eqnarray*}
Hence if $\left\| \frac{{\tilde A}}{N+2}\right  \| < 1$, then
\begin{eqnarray*}
 \left \| \sum_{m > N} A_{m} \right\| &\leq& \|A_N\|  \cdot
 \left\| \frac{{\tilde A}}{N+2}\right  \| \cdot
 \left(1 - \left\| \frac{{\tilde A}}{N+2}\right \|\right)^{-1}\\
& = & \|A_N\|  \cdot
 \frac{\|{\tilde A}\|}{N+2 - \| {\tilde A} \|} = r
\end{eqnarray*}

And finally,
\begin{equation}
 \sum_{m=0}^\infty A_m = \sum_{m=0}^N A_m + [-r, r]^n
\end{equation}

 {\bf An approach based on logarithmic norms:}(compare Lemma~\ref{lem:alg-lognorm})
  We fix
any norm $\|\cdot\|$, for example the $L^\infty$-norm:
$\|x\|_\infty=\max_i |x_i|$ (one should chose the norm which gives
the smallest $l$ )
\begin{itemize}
\item[1.]  $[\delta]=[\{f(x,y_c) - f(x,y) \: | \: x \in [W_1], y \in [W_y]  \}]_I$.
\item[2.] $C=\|[\delta]\|$
\item[3.] $l=\mbox{right}\left( \mu(\frac{\partial f}{\partial x}([W_2],y_c)) \right)$
\item[4.] If $l \neq 0$, then $D=\frac{C(e^{lh}-1)}{l}$. \newline
          If $l=0$, then $D=C h$
\item[5.] $[\Delta]=[-D,D]^n$
\end{itemize}

\noindent \textbf{ Remark.}
 In both cases we compute
\begin{equation}
 [\delta]=[\{f(x,y_c) - f(x,y) \: | \: x \in [W_1], y \in [W_y]
\}]_I. \label{eq:deltacomp}
\end{equation}
One need to be very careful in the computation  of $[\delta]$
using (\ref{eq:deltacomp}), because direct interval evaluation of
$[\{f(x,y_c) - f(x,y) \: | \: x \in [W_1], y \in [W_y] \}]_I$
yields big overestimation. Namely, when there is no perturbations
at all, i.e. $[W_y] = \{y_c\}$, then $[\delta]=0$. On the other
hand if $f([W_1]) = [\{f(x,y_c) \: | \: x \in [W_1]\}]_I =
[a^-,a^+]$ then  the naive interval computation give $[\delta] =
[a^- - a^+, a^+ - a^-]$, so $\diam [\delta] = 2  \diam f([W_1])$
and this can be big because $[W_1]$ is an enclosure of a solution
during the whole time step.

\subsection{Rearrangement}
The rearrangement is an essential ingredient in the Lohner
algorithm, designed to reduce the wrapping effect
\cite{Lo,Lo1,Mo}. We will not discuss this issue here, but we will
only include necessary formulas (see \cite{ZC1Lo} for more details
and the motivation).

\noindent \textbf{Evaluations 2 and 3.} In this representation
\begin{equation}
  [x_k]= x_k + [B_k][{\tilde r}_k].
\end{equation}
In the context of our algorithm in part 3 we obtain
\begin{equation}
  [\overline{x}_{k+1}] = \overline{x}_{k+1} + [B_{k+1}][{\overline r}_{k+1}].
\end{equation}
Now we have to take into account equation (\ref{eq:part4}). We set
\begin{eqnarray}
  x_{k+1}&=&{\rm m}(\overline{x}_{k+1} + [\Delta]) \label{eq:ev231} \\
  {}[\tilde r_{k+1}] &=& [\overline{r}_{k+1}] + [B_{k+1}^{-1}]\left(
    \overline{x}_{k+1} + [\Delta] - x_{k+1}  \right) \label{eq:ev232}.
\end{eqnarray}

 \noindent \textbf{Evaluation 4.} In this representation
 \begin{equation}
   [x_k]= x_k + C_k [r_0] + [B_k][{\tilde r}_k].
 \end{equation}
In the context of our algorithm in part 3 we obtain
\begin{equation}
  [\overline{x}_{k+1}] = \overline{x}_{k+1} + C_{k+1}[r_0] + [B_{k+1}][{\overline r}_{k+1}].
\end{equation}
Equation (\ref{eq:part4}) is taken into account exactly in the
same way as in previous evaluations, i.e., we use (\ref{eq:ev231})
and (\ref{eq:ev232}).

\section{Rigorous estimates between time steps}
\label{sec:betwper}

In order to compute the Poincar\'e map for differential inclusion
we also need an estimate for time $t \in [t_k, t_k + h_k]$.

 \vspace{0.5cm} \noindent \newline {\bf Input
parameters:}
\begin{itemize}
  \item $h_k$ is a time step
 \item $[x_k]\subset {\mathbb R}^n$, such that  $ \varphi(t_k,[x_0],[y_0]) \subset [x_k]$
 \item $[x_{k+1}]\subset {\mathbb R}^n$, such that  $\varphi(t_k+h_k,[x_0],[y_0]) \subset [x_{k+1}]$
 \item  convex and compact set $[W_2] \subset {\mathbb R}^n$  and  convex
      set  $[W_y] \subset {\mathbb R}^m$, such that
   \begin{eqnarray}
     \varphi([t_k,t_k+h_k],[x_0],[y_0]) \subset [W_2] \\
    {} [y_0]([t_k,t_{k+1}]) \subset [W_y].
  \end{eqnarray}
 \item $y_c \in [W_y]$
 \item $[\overline{x}_{k+1}] \subset {\mathbb R}^n$, such that
    $ \overline{\varphi}(h_k,[x_k],y_c) \subset  [\overline{x}_{k+1}]$
 \item $[W_1] \subset {\mathbb R}^n$ compact and convex, such that
    $\overline{\varphi}([0,h_k],[x_k],y_c) \subset [W_1] $
\end{itemize}

\noindent
 {\bf Output:} \newline
 We compute $[E_k] \subset {\mathbb R}^n$ such that
 \begin{displaymath}
    \varphi(t_k+[0,h_k],[x_0],[y_0]) \subset [E_k],
 \end{displaymath}

\noindent {\bf Algorithm:}
\begin{itemize}
 \item We compute $[\overline E_k] \subset {\mathbb R}^n$, such that
 \begin{equation}
    \overline{\varphi}([0,h_k],[x_k],y_c) \subset [\overline{E}_k]
 \end{equation}
 using a procedure for an ODE described in \cite{ZC1Lo}. This
 procedure requires as input data: $h_k$, $[x_k]$,
 $[\overline{x}_{k+1}]$ and $[W_1]$.

\item  we compute a set $[\Delta] \subset {\mathbb
R}^n$, such that
\begin{equation}
   \varphi(t_k + h,[x_0],[y_0]) \subset
  \overline{\varphi}(h,[x_k],y_c) + [\Delta], \qquad \mbox{for $0 \leq h \leq h_k$}.
\end{equation}
Observe that this requires $y_c$, $[W_1]$, $[W_2]$ and $[W_y]$.
\item finally we obtain
 \begin{equation}
     \varphi(t_k+[0,h_k],[x_0],[y_0])_i \subset [E_{k}]_i =
    [\overline{E}_{k}]_i + [\Delta]_i.
 \end{equation}
\end{itemize}

 \vspace{0.5cm}

Slightly better algorithm:
\begin{itemize}
 \item if $0 \notin f_i([W_2],[W_y])_i$, then the $i$-th coordinate is strictly
    monotone on $[W_2] \times [W_y]$, hence we  set
\begin{displaymath}
    [E_k]_i= \mbox{hull}([x_k]_i,[x_{k+1}]_i)
\end{displaymath}

 \item if $0 \in f_i([W_2],[W_y])$, then we compute $[\overline E_k] \subset {\mathbb
 R}^n$, such that
 \begin{equation}
   \overline{\varphi}([0,h_k],[x_k],y_c) \subset [\overline{E}_k]
 \end{equation}
 using a procedure for an ODE described in \cite{ZC1Lo}. This
 procedure requires as input data: $h_k$, $[x_k]$,
 $[\overline{x}_{k+1}]$and $[W_1]$.

 We have
 \begin{equation}
     \varphi(t_k+[0,h_k],[x_0],[y_0])_i \subset [E_{k}]_i =
    [\overline{E}_{k}]_i + [\Delta]_i.
 \end{equation}
\end{itemize}

\textbf{A drawback of this approach:} \newline if we have to
perform several time steps during which the computed enclosure for
the trajectory has a nonempty intersection with the section, then
$\Delta$ is added twice.

\subsection{Computation of the Poincar\'e map}
If as in \cite{ZC1Lo} we assume that the section is given by
$\alpha(x)=0$ then an algorithm discussed in Section 5 in
\cite{ZC1Lo} also applies  in the present context.

\section{Some tests, discussion}
\label{sec:testexamples}
\subsection{Perturbed harmonic oscillator}
\label{subsec:oscilator} We use  the harmonic oscillator to
compare two methods: first based on the logarithmic norms and the
second one that uses component-wise estimates. To shorten the
notation in this section we call them \emph{LN method} and
\emph{CW method} correspondingly.

The equations of the perturbed harmonic oscillator are given by
\begin{eqnarray}
 x' & = & y +\epsilon_1\label{eq:oscilator}\\
 y' & = & -x + \epsilon_2 \nonumber
\end{eqnarray}
and we will always use the initial condition given by $(1, 0)+
[-\delta, \delta]^2$.

In both methods we first find the solution of the unperturbed
system and then we add the influence of perturbation denoted
(following section \ref{sec:alg}) by $\Delta$. For this simple
system we are able to compute $\Delta$ for both methods by hand.
Let $h$ denote time step used.\\
For LN method we used the euclidian logarithmic norm $\mu_e$
because it is optimal for this case. Namely, we have
\begin{equation}
l = \mu_e(\frac{\partial f}{\partial x}([W_2],y_c) =
\mu_e\left( \left[\begin{array}{cc}
 0 & 1\\
 -1 & 0
\end{array}
 \right] \right)= 0.
\end{equation}
Therefore, we obtain $\Delta = [-D, D]^2$ where
\[
 D = h \sqrt{ \epsilon_1^2 +  \epsilon_2^2 }.
\]
For CW method we obtain $\Delta = ([-D_1, D_1], [-D_2, D_2])$,
where
\begin{eqnarray*}
 D_1 & = & \epsilon_1 \sinh h + \epsilon_2 (\cosh h - 1), \\
 D_2 & = & \epsilon_1 (\cosh h - 1) + \epsilon_2 \sinh h.
\end{eqnarray*}

Suppose that $\epsilon_1 = \epsilon_2 := \epsilon$, then LN method is better than CW method if
\begin{equation}
  \sqrt{2}h\epsilon < \epsilon(\sinh h + \cosh h -1) = \epsilon (\exp(h) -1)
\label{eq:lncw}
\end{equation}
Inequality (\ref{eq:lncw}) holds for $h > 0.657275$. As it can be
seen in Table~\ref{tab:osc1} results of computations agree with
this theoretical estimate and the  LN method is better for $h
> 0.657275$. We were  not able to use time steps
$h>0.8$ because for such a big time steps our rough enclosure
procedure (the first part of the algorithm) fails.

\begin{table}
\begin{center}
\begin{tabular}{|c|c|c|}
\hline
time step & LN method & CW method \\
$h$ & $D$ & $D_1, D_2$\\\hline
 0.799  &  0.112996  &  0.122332 \\
 0.7  &  0.0989949  &  0.101375 \\
 0.66  &  0.0933381  &  0.0934792 \\
 0.658  &  0.0930553  &  0.0930927 \\\hline
 0.657  &  0.0929138  &  0.0928997 \\
 0.65  &  0.0919239  &  0.0915541 \\
 0.5  &  0.0707107  &  0.0648721 \\
 0.25  &  0.0353553  &  0.0284025 \\
 0.1  &  0.0141421  &  0.0105171 \\
 0.01  &  0.00141421  &  0.00100502 \\
 0.001  &  0.000141421  &  0.00010005 \\\hline
\end{tabular}
\caption{ Perturbed harmonic oscillator $\epsilon_1=\epsilon_2=0.1$: Estimates of perturbations for various time steps - comparison between LN and CW method \label{tab:osc1}}
\end{center}
\end{table}

The situation is quite different, when we perturb only one
coordinate. Suppose that $\epsilon_1 = 0$ and $\epsilon_2 =
\epsilon$. Now, for LN method we have
\[
 D = h \epsilon,
\]
and for CW method
\begin{eqnarray*}
 D_1 & = &  \epsilon (\cosh h - 1) = \epsilon ( \frac{h^2}{2!} + \frac{h^4}{4!} + \dots), \\
 D_2 & = &  \epsilon \sinh h  = \epsilon ( h + \frac{h^3}{3!} + \frac{h^5}{5!} + \dots).
\end{eqnarray*}
From the above formulas it follows that for time steps up to
1.616137 value of $D_1$ is smaller than $D$, but $D_2$ is always
bigger than $D$. In Table~\ref{tab:osc2} we list values of
perturbations for LN an CW method for various time steps.  Again
for time steps bigger than 0.8 our implementation could not find
rough enclosure. For small time steps the ratio $\frac{D}{D_1}$ is
quite big, when the ratio $\frac{D}{D_2}$ is slightly less than
one. So overall it is better to use CW method.

\begin{table}
\begin{center}
\begin{tabular}{|c|c|c|c|}
\hline
time step & LN method & \multicolumn{2}{|c|}{CW method} \\
$h$ & $D$ & $D_1$ & $D_2$\\\hline
 0.8  &  0.08  &  0.0337435  &  0.0888106 \\
 0.5  &  0.05  &  0.0127626  &  0.0521095 \\
 0.25  & 0.025 &  0.0031413  &  0.0252612 \\
 0.1  &  0.01  &  0.0005004  &  0.0100167 \\
 0.01  &  0.001&  5.0e-06    &  0.0010001 \\
 0.001 & 0.0001&  5.002e-08  &  0.0001 \\\hline
\end{tabular}
\caption{ Perturbed harmonic oscillator $\epsilon_1=0, \epsilon_2=0.1$: Estimates of perturbations for various time steps - comparison between LN and CW method \label{tab:osc2}}
\end{center}
\end{table}

In Table \ref{tab:oscilator} we compare diameters of computed
rigorous estimates of solutions of (\ref{eq:oscilator}) after time
$T = 2 \pi$ for these two methods using various values of $h$,
$\epsilon$ and $\delta$. Again we perturb only second coordinate i.e. $\epsilon_1=0, \epsilon_2=\epsilon$.
 As expected, we see that decreasing time
steps results in the  increase of the accuracy of the estimates,
but it also increases computational cost. In the second part of
the table we were changing set sizes and in the third one we were
changing the size of the perturbation. It can be seen that our
algorithm is capable to provide estimates even for perturbations
much bigger than values of the vector field. Observe that
with the time steps used in these experiments the CW method is
better than LN method. The biggest time step $h$ used was approximately 
equal to $0.785$.

\begin{table}
\begin{center}
\begin{tabular}[c]{|r|r|r|c|c|}
\hline
& & number& \multicolumn{2}{c|}{size of the set after time $T = 2\pi$}\\
$\epsilon$ & $\delta$ & of steps & LN method & CW method \\ \hline
 0.1 & 0.01 & 8     &  1.5789308 & 1.2143687 \\
 0.1 & 0.01 & 100   &  1.6220657 & 0.8479880 \\
 0.1 & 0.01 & 1000  &  1.6202468 & 0.8227680\\
 0.1 & 0.01 & 10000 &  1.6200250 & 0.8202765\\
 0.1 & 0.01 & 100000&  1.6200025 & 0.8200276 \\ \hline

 0.1 & 0    & 100   &  1.5994735 & 0.8253958\\
 0.1 & 0.01 & 100   &  1.6220657 & 0.8479880\\
 0.1 & 0.1  & 100   &  1.8253953 & 1.0513176\\ \hline

 0.01& 0.01 & 100   &  0.1825395 & 0.1051317\\
 0.1 & 0.01 & 100   &  1.6220657 & 0.8479880 \\
 1   & 0.01 & 100   &  16.017328 & 8.2765505 \\
 10  & 0.01 & 100   &  159.96995 & 82.562176\\\hline
\end{tabular}
\caption{ Perturbed harmonic oscillator $\epsilon_1=0, \epsilon_2=\epsilon$: Estimates of perturbations for various values of the parameters - comparison between LN and CW method \label{tab:oscilator}}
\end{center}
\end{table}

\subsection{R\"ossler equations}
\label{subsec:Rossler}  R\"ossler equations \cite{R} are given by
\begin{eqnarray}
  x'&=&-(y+z)  \nonumber  \\
  y'&=&x+0.2y    \label{eq:rossler}  \\
  z'&=&0.2 + z(x-a),  \nonumber
\end{eqnarray}
where $a$ is a real parameter. In our tests  we set $a=5.7$ - the
'classical' parameter value for which numerical simulation display
a strange attractor \cite{R}.

In our test we focus on computation of a Poincar\'e map, $P$, on
section $\Theta=\{x=0, x'>0 \}$ around a point $x_0=(0.0, -10.3,
0.03)$. This is a point from the attractor (or close to the
attractor, which we have found numerically difficult in
\cite{ZHRo}).

In Table~\ref{tab:roosler} we list the results of a computation of
Poincar\'e map on section $\Theta$ for a differential inclusion
$x' \in f(x) + [\epsilon]$, where $f(x)$ is the vector field in
R\"ossler equations (\ref{eq:rossler}) and $[\epsilon]=[-10^{-4},
10^{-4}]^3$. The initial condition was $x_0 + \{0\} \times
[-10^{-4}, 10^{-4}]^2$. In computations the method based on the component-wise estimates and the Lohner algorithm - 4th evaluation was used. 

We see that our algorithm can provide good estimetes even for perturbed system and for set of initial data containing numerically difficult points from attractor.

\begin{table}
\begin{center}
\begin{tabular}{|l|c|}
\hline
initial set $[X]$& $(0.0, -10.3,
0.03) + \{0\}\times [-10^{-4},
10^{-4}]^2$\\
perturbations $[\epsilon]$ & $[-10^{-4},
10^{-4}]^3$\\ \hline
$P([X])$ &
 $\left(\begin{array}{c}
  [-0.3136278, 0.3049910]\\
 ~[-3.7425421,-3.4205722]\\
 ~[ 0.0306989, 0.0337781]
\end{array}\right)^T$\\
diam $P([X])$ & (0.6186189, 0.3219698,0.0030791) \\
 \hline
\end{tabular}
\caption{Perturbed R\"ossler equation: Value of a Poincar\'e map on
section $\Theta=\{x=0, x'>0 \} $\label{tab:roosler}}
\end{center}
\end{table}

\subsection{Kuramoto-Sivashinsky PDE's}
\label{subsec:KSodd} Assuming odd and periodic boundary conditions
the Kuramoto-Sivashinsky equations can be reduced \cite{ZM} to the
following infinite system of ordinary differential equations
\begin{equation}
\label{eq:KSfodd} \dot{a}_k  =  k^2(1-\nu k^2)a_k -
k\sum_{n=1}^{k-1} a_n a_{k-n}
 +2k\sum_{n=1}^\infty
a_n a_{n+k}\quad k=1,2,3,\ldots
\end{equation}
where $\nu >0$. In \cite{ZKSper,ZKS3} using the algorithm based on
component-wise estimates described in this paper to handle the
dominant modes and the method of self-consistent bounds developed
in \cite{ZM} to deal with the tail (the remaining modes) the
existence of multiple periodic orbits has been proved for a range
for $\nu \in [0.032,0.127]$. Some of these orbits were attracting,
while others were unstable with one unstable direction.


\end{document}